\documentclass[12pt]{article}
\input colordvi.sty
\usepackage{doi}
\usepackage[pdftex]{graphicx}
\usepackage{fancyhdr}
\usepackage{datetime}
\usepackage{cleveref}
\usepackage{amsfonts}
\usepackage{pstricks,pst-node,pst-text,pst-3d}
\usepackage{amsmath,amssymb,amsthm}
\usepackage{url}
\usepackage{hyperref}
\usepackage{amssymb}
\usepackage{color}

\usepackage{cleveref}
\usepackage{amsmath,amsfonts,amssymb,amscd,latexsym,lineno,epsfig}
\input colordvi.sty
\usepackage{graphicx}
\renewcommand{\ge}{\geqslant}

\newcommand{\N}{{\mathbb N}}

\newcommand{\be}{\begin{equation}}
\newcommand{\en}{\end{equation}}

\begin{document}
\numberwithin{equation}{section}

\begin{center}
\Large{{\bf KdV-type equations linked via B\"acklund transformations:  
remarks and perspectives
}}
\end{center}
\normalsize
\begin{center}{
{\bf\large Sandra Carillo}
\\{\rm Dipartimento Scienze di Base e Applicate
    per l'Ingegneria \\  \textsc{Sapienza}   Universit\`a di Roma,  16, Via A. Scarpa, 00161 Rome, Italy} \\  
    \& \\ {\rm National Institute for Nuclear Physics (INFN), Rome, Italy \\ Gr. Roma1, IV - Mathematical Methods in NonLinear Physics}\\ 
   \today}
 \end{center}
\medskip
\begin{abstract}  
Third order nonlinear evolution equations, that is the Korteweg-deVries (KdV),  modified  Korteweg-deVries (mKdV) equation  and o\-ther ones are considered: they all
are connected  via B\"acklund  transformations. These links can be depicted in a wide {\it B\"acklund Chart} which further extends the previous one constructed in \cite{Fuchssteiner:Carillo:1989a}. 
  In particular, the  B\"acklund transformation which links the mKdV equation to the  KdV singularity manifold equation is reconsidered and the nonlinear  equation for the KdV eigenfunction is shown to be linked to all the equations  in the previously constructed B\"acklund Chart.   That is, such a B\"acklund Chart is expanded to encompass  the nonlinear  equation for the KdV 
  eigenfunctions \cite{boris90}, which finds its origin in the early days of the study of   Inverse scattering Transform method, when the Lax pair for the KdV equation was constructed.   The nonlinear  equation for the KdV   eigenfunctions
    is proved  to enjoy a nontrivial invariance property. 
Furthermore, the hereditary recursion operator it admits  \cite{boris90}  is recovered via a different method. 
   Then, the  results are extended to the whole hierarchy of nonlinear evolution equations 
   it generates.
Notably, the established links allow to show that also the nonlinear  equation for the KdV   eigenfunction is connected 
to the  Dym equation since both such equations appear in the same B\"acklund  chart.
\end{abstract}

\noindent


\section{Introduction}
\label{introd}
The  relevance  of B\"acklund  transformations in soliton theory is well established, see  \cite{CalogeroDegasperis, RogersShadwick, RogersAmes, RogersSchief, Gu-book}
where a wide variety  of applications of B\"acklund and Darboux Transformations and their 
connections with partial differential equations admitting soliton solutions is given.
Here, the  concern is on  B\"acklund transformations as a tool to investigate structural 
properties of nonlinear evolution equations. Specifically, the B\"acklund  chart in 
 \cite{Fuchssteiner:Carillo:1989a}
is reconsidered to show that it can be further extended to incorporate also  the nonlinear  equation for the KdV eigenfunction. The latter,  
named hereafter {\it KdV eigenfunction} equation for short,
is studied in \cite{boris90} where, among many other ones,  it is proved to be integrable 
via  the {\it inverse spectral transform} (IST) method. Indeed, this equation was firstly 
derived  in a founding article of the IST method \cite{MGK} and also \cite{russi2}, later further investigated in  \cite{boris90, russi} wherein a  wide variety of nonlinear evolution equations is studied. Nevertheless, 
the KdV eigenfunction  equation does not appear in subsequent  classification studies of integrable nonlinear 
evolution equations, such as 
\cite{Calogero1985, Wang, Mikhailov-et-al, Mikhailov-et-al2} until very recently when, in \cite{Faruk1},  linearizable 
 nonlinear evolution equations are classified.
The KdV eigenfunction equation is a third order nonlinear equation of KdV-type since  it
is connected via B\"acklund  transformations with the  
Korteweg deVries (KdV),  the modified Korteweg deVries (mKdV), 
 the {\it Korteweg deVries interacting soliton} (int.sol.KdV)\cite{Fuchssteiner1987} and the 
 {\it Korteweg deVries singuarity manifold}  (KdV sing.)\footnote{The Korteweg deVries 
singuarity manifold equation, introduced in \cite{Weiss}, is also known as as UrKdV 
or Schwarz-KdV \cite{Depireux,  Wilson}.}  equations. The KdV eigenfunction equation is, then,
 proved to enjoy an  invariance property.  In addition,  since it is
connected via B\"acklund transformations to the other KdV-type  equations, according to 
 \cite{FokasFuchssteiner:1981, Fuchssteiner1979}, its  hereditary recursion operator \cite{boris90} 
 can be recovered.
The heritariness of all the recursion operators admitted by the equations in the B\"acklund  
chart allow to extend  all the links to 
the whole corresponding hierarchies; hence,  previous results  \cite{Fuchssteiner:Carillo:1989a}
are generalised. Generalisations  to non-Abelian KdV-type  equations and hierarchies are comprised in \cite{Carillo:Schiebold:JMP2009,  Carillo:Schiebold:JMP2011, 
SIGMA2016, JMP2018}, wherein the links among them are depicted in a noncommutative B\"acklund  
chart analogous of that one in  \cite{Fuchssteiner:Carillo:1989a}.

\medskip
The material is organized as follows.
The opening Section \ref{background} is devoted to recall the definition of B\"acklund  transformation 
adopted throughout this work together with its consequences which are most relevant  to the present 
investigation. In the following 
Section \ref{new-eq} the nonlinear equation for the KdV eigenfunction is obtained. Specifically, it is 
shown to 
be linked, via  B\"acklund  transformations with the mKdV and the KdV singularity manifold equations. 
Notably, both the equations introduced in \cite{JMP2018}  represent non-Abelian counterparts  of this equation when commutativity is assumed.
In Section \ref{inv},   the KdV eigenfunction equation is proved to enjoy a 
non trivial invariance property. The subsequent Section  \ref{ext-bc} concerns   
the  B\"acklund  chart in \cite{Fuchssteiner:Carillo:1989a}, which is  further extended to include also
the  KdV eigenfunction equation. 
In Section \ref{rec-op}, via the links in the B\"acklund  chart, the   
hereditary recursion operator, firstly obtained in \cite{boris90}, admitted by the KdV eigenfunction 
equation is recovered in explicit form. 
Thus,  the generated hierarchy follows.
As a consequence \cite{FokasFuchssteiner:1981}, the  whole hierarchy of nonlinear evolution 
equations turns out to be connected to all the hierarchies in  the B\"acklund  chart. 
Concluding remark  as well as how the present work is related to previous ones and, in particular, with the research program devoted to the study of non-Abelian nonlinear evolution equations \cite{Carillo:Schiebold:JMP2009, Carillo:Schiebold:JMP2011, SIGMA2016, JMP2018},  are comprised in the closing Section \ref{rem}.
\section{Some background definitions}
\label{background}
This section is devoted to provide some background definitions which are of use throughout the whole 
article. Since many definitions are not unique in the literature, here, those ones here adopted are provided. 

First of all,  the  notion of B\"acklund Transformation, according to Fokas and Fuchssteiner  
\cite{FokasFuchssteiner:1981} is recalled (see also the book by Rogers and Shadwick \cite{RogersShadwick}). 
Consider non linear evolution equations of the type
\begin{equation}\label{1}
u_t = K ( u ) 
\end{equation}
where the unknown function $u$ depends on the independent variables $x,t$ and, for fixed 
$t$, $u (x,t) \in M$, a manifold modeled on a linear topological space so that the generic 
{\it fiber} $T_uM$, at $u\in M$, can be identified 
with $M$ itself \footnote{It is generally assumed that $M$ is the space of functions $u(x,t)$ which, 
for each fixed $t$, belong to the Schwartz space $S$ of {\it rapidly decreasing functions} on 
${{\mbox{R\hskip-0.9em{}I \ }}}^n$,    i.e.
$S({{\mbox{R\hskip-0.9em{}I \ }}}^n):=\{ f\in C^\infty({{\mbox{R\hskip-0.9em{}I \ }}}^n) : \vert\!\vert f \vert\!
\vert_{\alpha,\beta} < \infty, \forall \alpha,\beta\in \N_0^n\}$, where 
$\vert\!\vert f \vert\!\vert_{\alpha,\beta}:= sup_{x\in{{\mbox{R\hskip-0.9em{}I \ }}}^n} \left\vert x^\alpha D^\beta f(x)
\right\vert $, and  $D^\beta:=\partial^\beta /{\partial x}^\beta$; throughout this article $n=1$.}.  and  $ K : M 
\rightarrow TM$,  is an appropriate $ C^{\infty}$ vector field on a manifold $ M$ to its  
tangent manifold $TM$.  Let, now  
\begin{equation}\label{2}
v_t = G (v)
\end{equation}
denote another nonlinear evolution equation. If it assumed $ u (x,t) \in M_1$ and $ v (x,t) \in M_2 $ 
where $M_1, M_2$ represent  manifolds modeled on a linear topological space 
so that, $ K : M_1 \rightarrow TM_1 $ and $ G : M_2 \rightarrow TM_2 $ represent  appropriate 
$ C^{\infty}$-vector fields on the manifolds $M_i, i=1,2$, then \vskip-1em
\begin{eqnarray}\label{eq.s}
u_t &= K ( u ),~~~ K : M_1 \rightarrow TM_1,~~~{{ u  :(x,t)    \in{{\mathbb R}} 
 \times  {{\mathbb R}}\to u (x,t)   \in M_1}}\\
v_t &= G (v),~~~G : M_2  \rightarrow TM_2,~~~v  :(x,t)    \in{{\mathbb R}} 
 \times{{\mathbb R}}  \to   v (x,t)   \in  M_2 .
\end{eqnarray}
Here, according to the usual choice when soliton solutions are considered, it is further assumed 
$M:= M_1\equiv M_2$.
Then,  \cite{FokasFuchssteiner:1981}\,   a B\"acklund transformation can defined as follows.

\medskip\noindent
{\bf{ Definition}} {\it Given two evolution equations, $ u_t = K (u)$ and $v_t = G (v)$, then $\hbox{B (u , v) = 0}$  
represents a B\"acklund transformation between them 
 whenever  given two solutions of such  equations, say, respectively, $u(x,t)$ and $v(x,t)$  such that 
\begin{equation}
B (u(x,t), v(x,t)) \vert_{ t=0 } = 0 
\end{equation}
it follows that,
\begin{equation}
B (u(x,t),v(x,t) )\vert_{t=\bar t} = 0,     ~~\forall \bar t >0 ~, ~~~\forall x\in{\mathbb R}.
\end{equation}
}
\noindent Hence, solutions admitted by the two equations are connected via the B\"acklund transformation 
which establishes a correspondence between them: it can  graphicallly represented as %
\begin{eqnarray}\label{BC1}
\boxed{u_t = K (u)} \,{\buildrel B \over {\textendash\textendash}}\, \boxed{ v_t = G (v) }~~.
\end{eqnarray}
In addition, if,  the nonlinear evolution equation  (\ref{1}) admits a hereditary recursion operator 
\cite{FokasFuchssteiner:1981, Olver}, \,,   denoted as $\Phi (u)$,  it  can be written as
\begin{equation}\label{base_u}
u_t = \Phi( u ) u_x ~~~\text{where}~~~ K (u) = \Phi (u) u_x~.
\end{equation}
The B\"acklund  transformation allows to construct the operator  $\Psi(v)$ 
\begin{equation}\label{transf-op} 
\Psi(v)=  \Pi\, \Phi (u)\, \Pi^{ -1}  
\end{equation}
where
\begin{equation}\label{pi}
 \Pi : = -B_v^{ -1} B_u~~,~~
          \Pi : T M_1  \rightarrow T M_2~,
\end{equation}
while $B_u$ and $B_v$ denote the Frechet derivatives of the B\"acklund transformation $B(u,v)$.
Then, \cite{FokasFuchssteiner:1981}\,, the operator $\Psi(v)$ represents the  hereditary recursion operator admitted by   the   equation $v_t=G(v)$ which, thus, can be written under the form
$$G(v) = \Psi (v)\, v_x~.$$ 
That is, according to \cite{FokasFuchssteiner:1981}, given the B\"acklund transformation $B(u,v)$, 
and the hereditary recursion operator  $\Phi (u)$ admitted by equation \eqref{1}, then, also  
equation \eqref{2} admits a hereditary recursion operator: it is obtained on
 use of the operator $\Pi$, \eqref{pi}, via the trasformation formula \eqref{transf-op}.
 
On subsequent applications of the admitted recursion operators are, 
respectively, the  two hierarchies
\begin{equation}
\displaystyle u_{t} =\left[\Phi (u)\right]^{n} u_x ~~~~~\text{and} ~~~~~~ 
v_{t }= \left[\Psi (v)\right]^{n} v_x~,~~ n\in\N 
\end{equation}
of  evolution equations can be constructed   \cite{Fuchssteiner1979}; their {\it base members} 
equations, which correspond to $n=1$, coincide with  equations (\ref{1}) 
and (\ref{2}). Fixed any $n_0\in\N$, the two equations $u_{t} =\left[\Phi (u)\right]^{n_0} u_x$ and $v_{t }= \left[\Psi (v)\right]^{n_0}$ are connected, via the same 
B\"acklund Transformation  which connects the two base members equations.
This extension to the whole hierarchies is graphically represented by the following 
B\"acklund chart  
\begin{equation}\label{BT-hier}
\boxed{u_{} = \left[\Phi (u)\right]^{n} u_x }{\,\buildrel B \over {\textendash\textendash}}\,
\boxed{ v_{t} \ =\  \left[\Psi (v)\right]^{n} v_x} ~.   
\end{equation}
which emphasizes that 
the link between the two equations (\ref{1}) and (\ref{2}) is  extended  to corresponding members 
of  the two hierarchies generated, respectively, by the recursion operators $\Phi$ and $\Psi$.

\section{A third order KdV-type equation}
\label{new-eq}
In this Section the B\"acklund  chart  in \cite{Fuchssteiner:Carillo:1989a} 
is further extended to incorporate the  KdV eigenfunction equation:
 \begin{equation}\label{new}
w_t = w_{xxx} - 3 {{w_x w_{xx}}\over w} ~.
\end{equation}
This non linear evolution equation was introduced in  \cite{MGK}, one of the IST founding articles.  Later, it was
investigated in  \cite{russi2} and, subsequently,  in \cite{boris90} \footnote{see also  \cite{BSKonop}}, where its 
integrability via the {\it inverse spectral transform} (IST) method of soliton eigenfunction equations is proved.  Among 
the many equations studied in the extended article \cite{boris90} also the  KdV eigenfunction equation is included:  
ist  IST integrability is proved and also, via Lax pair representation, its recursion operator is provided. 
In Section \ref{rec-op} the explicit form of the recursion operator admitted by (\ref{new}) is constructed 
via a different approach. Indeed, to obtain such a recursion operator, we apply the connections, via 
B\"acklund transformations, of equation (\ref{new}) with other KdV-type equations according to the result presented in Section \ref{ext-bc} and here. Notably,  equation (\ref{new})  appears also in recent works \cite{russi2, Faruk1}. The latter
finds this equation in classifying linearizable evolution equations.

In this Section,  equation  \ref{new} is shown to be linked with the mKdV and the  KdV  singularity manifold equations. The following proposition can be proved. \medskip

\noindent {\bf Proposition 1} \\ \noindent 
{\it Equation (\ref{new}) is linked to the mKdV equation
\begin{equation}\label{mkdv}
v_t = v_{xxx} - 6 v^2  v_x 
\end{equation}
via the Cole-Hopf {\rm\cite{Cole:1951, Hopf:1950}} transformation
\begin{equation}\label{CH}
\text{\rm CH}: ~~~~vw- { w_x}=0~.
\end{equation}
} 
\noindent {\bf Proof} \\ \noindent 
  On substitution,  in (\ref{mkdv}), of $v$ in terms of $w$ according to the latter gives:
\begin{equation*}
\left( { w_x\over w}\right)_t = \left({ w_x\over w}\right)_{xxx} - 6 \left({ w_x\over w}\right)^2 \left({ w_x\over w}\right)_x
\end{equation*}
which, since the assumed regularity of $v$ and $w$ implies  Schwartz theorem on order of partial derivation holds, delivers 
\begin{equation*}
\left( { w_t\over w}\right)_x = \left[{ w_{xxx}\over w}- 3{{ w_xw_{xx}}\over w^2} +2  \left({ w_x\over w}\right)^3 -2  \left({ w_x\over w}\right)^3\right]_x
\end{equation*}
hence, after simplification, equation (\ref{new}) follows.  \hfill$\Box$
\medskip

\noindent {\bf Proposition 2} \\ \noindent 
{\it Equation (\ref{new}) is linked to the KdV sing manifold equation
\begin{equation}\label{phi-eq}
\varphi_t =  \varphi_x  \{ \varphi ; x\} ~~,~~ \{ \varphi ; x \} = \left( { \varphi_{xx} \over \varphi_x} \right)_x - {1 \over 2 }\left({ \varphi_{xx} \over \varphi_x} \right)^2
\end{equation}
via the B\"acklund transfornation
\begin{equation}\label{hatB}
\text{\rm B}: ~~~~\displaystyle{w^2 -  { \varphi_x}=0~.}
\end{equation}
} 
\noindent {\bf Proof} \\ \noindent 
The B\"acklund  transformation \eqref{hatB} implies:
\begin{eqnarray*}
\displaystyle 2 w w_t &=& \varphi_{xt} \\
\displaystyle 2{w_{xx} \over w} &=&  {\varphi_{xxx} \over \varphi_x} -{1\over 2} {\varphi_{xx}^2 \over \varphi_x^2}\\
\displaystyle 2{w_{xxx} \over w} &=&  {\varphi_{xxxx} \over \varphi_x} -{3\over 2} {{\varphi_{xx} \varphi_{xxx}} \over \varphi_x^2} 
+{3\over 4} {\varphi_{xx}^3 \over \varphi_x^3}
\end{eqnarray*}
Substitution of the latter in (\ref{new}) gives:
\begin{eqnarray*}
\displaystyle \varphi_{xt} = \varphi_{xxxx}  -{3} {{\varphi_{xx} \varphi_{xxx}} \over \varphi_x}  +{3\over 2} {\varphi_{xx}^3 \over \varphi_x^2}
\end{eqnarray*}
since \begin{equation*}
\displaystyle  \{ \varphi ; x \}_x = \varphi_{xxxx}  -{3} {{\varphi_{xx} \varphi_{xxx}} \over \varphi_x}  +{3\over 2} {\varphi_{xx}^3 \over \varphi_x^2}
\end{equation*}
on integration with respect to $x$, the KdV  sing manifold equation \eqref{phi-eq}  follows and the proof is complete.  \hfill$\Box$
\smallskip

\medskip

\noindent {\bf Remark} \\ 
\noindent 
Both the two new non-Abelian nonlinear evolution equations\footnote{Capital case is used to emphasize that the
unknown fuctions $Z$ and $W$ are non-Abelian ones, in \cite{JMP2018} operators on a Banach space.}
\begin{equation}\label{ckdv}
W_t = W_{xxx} - 3\, W_{xx}\, W^{-1} W_x ~, ~~ Z_t = Z_{xxx} - 3 \,Z_x \, Z^{-1} Z_{xx}~ 
\end{equation}
 introduced in \cite{JMP2018}, reduce to equation \eqref{new} 
when the assumption of a non-Abelian unknown is  removed. 

 \section {A non trivial invariance}
 \label{inv} 
 
Some properties  the introduced  equation \eqref{new} enjoys are studied in this section and in the 
following ones. In particular,  this section is devoted to prove an invariance property it enjoys. 
 
 First of all,  equation \eqref{new} is {\it scaling invariant} since 
  substitution of  $\alpha w, \forall \alpha \in \mathbb{C}$, to $w$ leaves it unchanged.  
  In addition, the following  the following proposition holds.
\medskip

\noindent {\bf Proposition 2} \\ \noindent 
{\it  The nonlinear evolution equation \eqref{new} is  invariant under the transformation 
\begin{equation}
\text{\rm I}: ~~~ \hat w^2 ={{ad- bc}\over{(c D^{-1}( w^2) +d)^2}}w^2,\quad a,b,c,d\in \mathbb{C}   
~ \text{s.t.} ~ad-bc\neq 0, 
\end{equation}
}
where
\begin{equation*}
D^{-1}:=\int_{-\infty}^x d\xi 
\end{equation*}
is well defined since so called {\it soliton solutions}  are looked for in the  Schwartz space $S({\mathbb R}^n)$ \footnote{see footnote on page 4.}.

\noindent {\bf Proof} \\ \noindent 
The KdV singularity manifold equation \eqref{phi-eq}  is  invariant under the M\"obius group of 
transformations
\begin{equation}
\text{M}:~~ \hat\varphi={{a\varphi+b}\over{c\varphi+d}},\qquad a,b,c,d\in \mathbb{C} \qquad \text{such that} \quad ad-bc\neq 0.
\end{equation}
Recalling that, according to proposition 2, the KdV singularity manifold 
equation \eqref{phi-eq} and equation \eqref{new} are connected to each other  via the B\"acklund  transformation B \eqref{hatB}, the result is readily obtained. 
Indeed,  the following B\"acklund chart 
 \begin{eqnarray*}
\text{M}: \hat \varphi={{a\varphi+b}\over{c\varphi+d}}~~~~~~~~~~~~~~~~\boxed{\!w_t = w_{xxx} - 3 {{w_x w_{xx}}\over w}\! } ~\, {\buildrel {w^2 -  { \varphi_x}=0} \over{\text{\textendash\textendash\textendash\textendash\textendash\textendash}}}~\boxed{\varphi_t \ =\  \varphi_x  \{ \varphi ; x\}}\\
\updownarrow~ \text{I} ~~~\qquad\qquad\qquad\qquad~~~~~\updownarrow~ \text{M} ~~~~~ \\
\forall a,b,c,d\in \mathbb{C} \vert ~ad-bc\ne 0~~~~~ \boxed{\!\hat w_t = \hat w_{xxx} - 3 {{\hat w_x \hat w_{xx}}\over \hat w}\! } \,~ {\buildrel {\hat w^2 -  { \hat\varphi_x}=0 }\over{\text{\textendash\textendash\textendash\textendash\textendash\textendash}}}~\boxed{\hat\varphi_t \ =\  \hat \varphi_x  \{ \hat \varphi ; x\}}\\
\end{eqnarray*}
shows that the invariance I  is obtained via composition of the M\"obius   transformation M with the two B\"acklund transformations
\begin{equation*}
\text{\rm B}: ~~~~\displaystyle{w^2 -  { \varphi_x}=0~ }~~~~\text{and}~~~~~~~~~~\widehat {\text{\rm B}}: ~~~~ \displaystyle{\hat w^2 -  {\hat \varphi_x}=0~.}
~~~ \qquad\qquad\qquad\Box
\end{equation*}
\section {Extended B\"acklund chart}
 \label{ext-bc}
 In this section, the equation \eqref{new}  is inserted in the B\"acklund chart in 
 \cite{Fuchssteiner:Carillo:1989a} which, then, is further extended.  
 Indeed, combination of the two transformations CH \eqref{CH} and B 
 in the previous section gives
 \begin{equation}\label{vphi}
\displaystyle v -{ {1 \over 2 }{ \varphi_{xx} \over \varphi_x} } = 0
\end{equation}
which coincides with the link  in  \cite{Fuchssteiner:Carillo:1989a} between mKdV and KdV sing. equations\footnote{Combination of (1.12) with (1.19), respectively, Cole-Hopf and introduction of a {\it bona fide} potential when connecting the  interacting soliton KdV with the  KdV sing. equation  \cite{Fuchssteiner:Carillo:1989a} produce the transformation \eqref{vphi}.}. 
The links given up to this stage  are sumarised in the following  B\"acklund chart: 
\begin{gather*}
\mbox{$\footnotesize
\boxed{\!\text{KdV}(u)\!}\, {\buildrel (a)
\over{\text{\textendash\textendash}}} \, \boxed{\!\text{mKdV}(v)\! } \, {\buildrel (b) \over{\text{\textendash\textendash}}} \, \boxed{\!\text{new eq}(w)\! } \, {\buildrel (c) \over{\text{\textendash\textendash}}} \, \boxed{\!\text{KdV~sing.}(\varphi)\!}
 {\buildrel (d) \over{\text{\textendash\textendash}}}\, \boxed{\!\text{int. sol KdV}(s) \!} \,
{\buildrel (e) \over{\text{\textendash\textendash}}} \, \boxed{\!\text{Dym}(\rho)\!}\,$}\label{BC1*}
\end{gather*}
where all the third order nonlinear evolution equations are, respectively 
\begin{alignat*}{3} 
& u_t = u_{xxx} + 6 uu_x \qquad && \text{(KdV)}, & \\
& v_t = v_{xxx} - 6 v^2 v_x \qquad && \text{(mKdV)},& \\
& w_t = w_{xxx} - 3 {{w_x w_{xx}}\over w}\qquad && \text{(new eq)},& \\
& \varphi_t = \varphi_x \{ \varphi ; x\} , \quad \text{where} \ \ \{ \varphi ; x \} :=
 \left( { \varphi_{xx} \over \varphi_x} \right)_x -
{1 \over 2 }\left({ \varphi_{xx} \over \varphi_x} \right)^2 \qquad && \text{(KdV~sing.)}, & \\
& s^2 s_t = s^2 s_{xxx} - 3 s s_x s_{xx}+ {3 \over 2 }{s_x}^3 \qquad && \text{(int.\ sol~KdV)}, & \\
& \rho_t = \rho^{3} \rho_{\xi \xi \xi} \qquad && \text{(Dym)}.&
\end{alignat*}
The B\"acklund transformations  linking these equations them, are, following their order in the B\"acklund chart:
 \begin{gather}\label{links}
(a) \ \ u + v_x + v^2 =0 , \qquad \qquad \qquad (b) \ \ v -{{ w_{x} \over w} } = 0,\\
(c) \  w^2 -{  \varphi_x}  = 0,\!\qquad \qquad  \qquad\qquad (d) \ \ s - \varphi_x =0, 
\end{gather} 
and 
\begin{equation}
(e) ~\ {\bar x} : = D^{-1} s (x), ~\rho(\bar x) :=  s(x), ~~~~~\text{where}~~~ ~D^{-1}:= \int_{-\infty}^x d\xi , \label{rec}
\end{equation}
so that  $\bar x= \bar x(s,x)$ and, hence, $ \rho(\bar x) :=  \rho(\bar x(s,x))$. The transformation  
(e) denotes the reciprocal transformation, as it is termed to stress it interchanges the role of the 
dependent and independent variables\footnote{ see, for instance,  \cite{RogersShadwick} for a 
general introduction and application of   reciprocal transformations. Details on the transformation $
(e)$ are given in \cite{BS1, Fuchssteiner:Carillo:1989a}. }
It represents a B\"acklund transformation between 
the extended manifold consisting of the both the dependent and the independent variables;
namely, the manifold given by pairs
formed by the dependent and the independent variables.
Then, the transformation (e) defines, at least locally \cite{BS1, FokasFuchssteiner:1981}
a B\"acklund transformation:
 $$ T_{(s,x)}   \rightarrow   T_{ ( \rho , \bar x ) }    $$
between the two respective tangent spaces. Therefore, it
is possible to transfer
the infinitesimal structure using the transformation formulae
 \cite{Fuchssteiner:Carillo:1989a, FokasFuchssteiner:1981}.  
 
Now, taking into account the invariance under the M\"obius group of transformations enjoyed by the 
singularity manifold equation, the B\"acklund chart can be  duplicated to obtain
\begin{eqnarray*}
\mbox{\footnotesize $ \boxed{\!\text{KdV}(u)\!}\, {\buildrel (a)
\over{\text{\textendash\textendash}}} \, \boxed{\!\text{mKdV}(v)\! } \, {\buildrel (b) \over{\text{\textendash\textendash}}} \, \boxed{\!\text{new eq}(w)\! } \, {\buildrel (c) \over{\text{\textendash\textendash}}} \, \boxed{\!\text{KdV~sing.}(\varphi)\!}
 {\buildrel (d) \over{\text{\textendash\textendash}}}\, \boxed{\!\text{int. sol KdV}(s) \!} \,
{\buildrel (e) \over{\text{\textendash\textendash}}} \, \boxed{\!\text{Dym}(\rho)\!}\,$}
 \\ \footnotesize 
AB_1\updownarrow~ ~~~ ~~AB_2 \updownarrow~~~~~~~~AB_3 \updownarrow~~~~~~~~~~~M \updownarrow~~~~~~~~~~~~AB_4 \updownarrow~~~~~~~~~~~~AB_5 \updownarrow~~~ \\
\mbox{\footnotesize $\boxed{\!\text{KdV}(u)\!}\, {\buildrel (a)
\over{\text{\textendash\textendash}}} \, \boxed{\!\text{mKdV}(v)\! } \, {\buildrel (b) \over{\text{\textendash\textendash}}} \, \boxed{\!\text{new eq}(w)\! } \, {\buildrel (c) \over{\text{\textendash\textendash}}} \, \boxed{\!\text{KdV~sing.}(\varphi)\!}
 {\buildrel (d) \over{\text{\textendash\textendash}}}\, \boxed{\!\text{int. sol KdV}(s) \!} \,
{\buildrel (e) \over{\text{\textendash\textendash}}} \, \boxed{\!\text{Dym}(\rho)\!}\,$}
\end{eqnarray*}
where the vertical lines represent auto-B\"acklund transformations which are all induced by the  
combination of  the B\"acklund transformations linking the other equations with invariance enjoyed 
by the KdV~singularity equation.  In detail, starting from the left hand side, the invariance 
AB$_1$ and AB$_2$ are the well known KdV and mKdV \cite{Miura},  auto-B\"acklund transformations, 
given in \cite{CalogeroDegasperis, Fuchssteiner:Carillo:1989a}, AB$_3\equiv$ I  is the invariance 
admitted by the novel  nonlinear evolution equation \eqref{new}  proved in Proposition 2.   
The last two, AB$_4$, and  AB$_5$, respectively, are auto-B\"acklund transformations,  
\cite{Fuchssteiner:Carillo:1989a}, admitted by the int. sol. KdV and Dym equations.
Notably, the connection between KdV and Dym equation \cite{RogersNucci, Fuchssteiner:Carillo:1989a} 
finds applications in the construction of solutions admitted by the Dym equation  \cite{BS4, Guo:Rogers}. 
The  extension to $2+1$ dimensional equations is given in  \cite{Rogers:1987, walsan1}.

\section {Admitted recursion operator \& hierarchy}
 \label{rec-op}

This section is devoted to the construction of the recursion operator admitted by the equation \eqref{new}. Specifically, the 
Cole-Hopf link  \eqref{hatB} between equation \eqref{new} and the mKdV equation allows to prove, according to 
\cite{Fuchssteiner1979, FokasFuchssteiner:1981}, it admits a recursion operator.

\medskip

\noindent {\bf Proposition 3} \\ \noindent 
{\it  The nonlinear evolution equation \eqref{new} admis  the recursion operator    
\begin{equation}\label{psiw}
\Psi(w)= {1 \over {2 w}}D w^2 \left[ D^2+2 U + D^{ -1} 2 UD\right] {1 \over {w^2 }} D^{ -1} 2w~~,
\end{equation}
where 
\begin{equation}\label{U}
U:=  {w_{xx} \over {w}}- 2 {w_{x}^2 \over {w^2}}~~.
\end{equation}
}
\noindent {\bf Proof} \\ \noindent 
Consider the B\"acklund transformation \eqref{hatB}, then, the transformation operator $\Pi$, recalling \eqref{pi}, reads 
\begin{equation}
\Pi : = -B_w^{ -1} B_{\varphi}
\end{equation}
where 
\begin{equation}
\displaystyle B_w[q]=\left. {\partial \over {\partial\varepsilon}} B(w+\varepsilon q, \varphi) \right \vert_{\varepsilon=0}=\left.
{\partial \over {\partial\varepsilon}}\left[ (w+\varepsilon q)^2- \varphi_x\right] 
\right\vert_{\varepsilon=0}=  2wq 
\end{equation}
and, 
\begin{equation}
\displaystyle B_{\varphi}[q]=\left. {\partial \over {\partial\varepsilon}} B(w, {\varphi}+\varepsilon q) \right \vert_{\varepsilon=0}=
\left.{\partial \over {\partial\varepsilon}}\left[ w^2+( \varphi+\varepsilon q)_x\right] 
\right \vert_{\varepsilon=0}=  q_x
\end{equation}
hence
\begin{equation}
B_w= 2 w~~~, ~ B_{\varphi} = D~\Rightarrow~\Pi= -{B_w}^{-1}B_{\varphi} = {1 \over{2w}}D~.
\end{equation}
Now, substitution of the latter in \eqref{transf-op} gives
\begin{equation}\label{12}
\displaystyle \Psi(w)= \left.{B_w}^{-1}B_{\varphi} \Phi(\varphi) B_{\varphi}^{-1}{B_w}\right\vert_{w^2 -{  \varphi_x}  = 0},\end{equation}
where  $\Phi(\varphi)$,   according to formulae (1.23)-(1.24) in \cite{Fuchssteiner:Carillo:1989a}, is 
 \begin{equation}
\displaystyle  \Phi({\varphi})=  \varphi_x\left[ D^2 +  \{ \varphi ; x\}+D^{-1}  \{ \varphi ; x\}\right] {1\over \varphi_x} D^{-1}~~.
\end{equation}
Then, substitution of the latter and of the transformation operator $\Pi$ in \eqref{12}
gives the operator
\begin{equation}
\displaystyle \Psi({w})= \left.{1 \over{2w}}D \varphi_x\left[ D^2 +  \{ \varphi ; x\}+D^{-1}  \{ \varphi ; x\}D\right] {1\over \varphi_x} D^{-1}{2w}\right\vert_{\varphi_x=w^2}~~.
\end{equation}
recalling the definition \eqref{phi-eq}$_2$ of the Schwarzian derivative $\{ \varphi ; x\}$,  
on substitution of $\varphi_x=w^2$, it follows
\begin{equation}
\displaystyle \left. \{ \varphi ; x\} \right\vert_{\varphi_x=w^2}=  
2\left( {w_{xx} \over {w}}- 2 {w_{x}^2 \over {w^2}}\right), 
\end{equation}
the latter, on use of  $U$, introduced in \eqref{U} to simplify the notation, gives
\begin{equation*}
\displaystyle \left. \{ \varphi ; x\} \right\vert_{\varphi_x=w^2}=  
2  U
\end{equation*}
and, hence,
\begin{equation*}
\displaystyle \Psi(w)= {1 \over{2w}}D w^2\left[ D^2 +   2U+  D^{-1}  2U D\right] {1\over w^2} D^{-1}{2w}~~,
\end{equation*}
which coincides with \eqref{psiw} and completes the proof. $\qquad\qquad\qquad\qquad\quad\Box$

\medskip \noindent 
In addition, the following proposition holds. 
\noindent {\bf Proposition 5} \\ \noindent 
{\it  The recursion operator \eqref{psiw} admitted by the nonlinear evolution equation \eqref{new} is  hereditary.}   

\bigskip 
\noindent {\bf Proof} \\ \noindent 
To prove the thesis, note that equation \eqref{new} is linked   via B\"acklund transformations to all the nonlinear evolution equations in the B\"acklund chart; hence, since all of them admit a hereditary recursion operator, according to \cite{FokasFuchssteiner:1981, Fuchssteiner1979}, also the  recursion operator \eqref{psiw}  admitted by  the newly obtained equation  \eqref{new} enjoys the  hereditariness property.   \hfill$\Box$

\medskip
\noindent {\bf Remark} \\ \noindent 
{  The recursion operator \eqref{psiw} admitted by the equation \eqref{new}  can be also obtained 
from the Cole-Hopf link ${vw -{ w_x}  = 0}$, with the mKdV equation, via the same method.
Hence, in this case  
\begin{equation}
\displaystyle \Psi(w)= \left.{B_w}^{-1}B_{v} \Phi_{\text{mKdV}}(v) B_{v}^{-1}{B_w}\right\vert_{vw -{ w_x}  = 0},\end{equation}
where, respectively,  $B_w$ and $B_v$ are given by
\begin{equation}
\displaystyle B_w[q]=\left . {\partial \over {\partial\varepsilon}} B(w+\varepsilon q, \varphi) \right \vert_{\varepsilon=0}=\left.
{\partial \over {\partial\varepsilon}}\left[ (w+\varepsilon q)^2- \varphi_x\right] \right \vert_{\varepsilon=0}=  2wq =
\end{equation}
and, 
\begin{equation}
\displaystyle B_{v}[q]=\left . {\partial \over {\partial\varepsilon}} B(w, v+\varepsilon q) \right \vert_{\varepsilon=0}=\left.
{\partial \over {\partial\varepsilon}}\left[ w(v+\varepsilon q)-w_x\right] \right \vert_{\varepsilon=0}=  wq 
\end{equation}
 and
 \begin{equation}
\displaystyle B_{w}[q]=\left . {\partial \over {\partial\varepsilon}} B(w+\varepsilon q, v) \right \vert_{\varepsilon=0}=\left.
{\partial \over {\partial\varepsilon}}\left[ (w+\varepsilon q) v-(w+\varepsilon q)_x\right] \right \vert_{\varepsilon=0}=
  vq-q_x 
\end{equation}
Explicit computations allow to obtain, once again,    \eqref{psiw}.
}

\bigskip\noindent
Now, equation \eqref{new} , when  the hereditary recursion operator $\Psi(w)$ is given in  \eqref{psiw}, 
can be written  as
\begin{equation}
\displaystyle w_t = \Psi(w)  w_x~~~ 
\end{equation}
and the corresponding  hierarchy is generated
\begin{equation}\label{newh}
\displaystyle w_t = \left[\Psi(w) \right]^n w_x~~~, ~n\in\N.
\end{equation}
Since,   as in  \cite{Fuchssteiner:Carillo:1989a}, all the nonlinear evolution equations in  the  
B\"acklund chart in Section \ref{ext-bc} admit a hereditary recursion operator    
\cite{FokasFuchssteiner:1981, Fuchssteiner1979},   all  the links can be extended 
to the corresponding whole hierarchies. Then, fixed $n=n_0,  n_0\in\N$ a different B\"acklund chart is obtained which links nonlinear evolution equations of order $2n_0+1$;
the case  $n_0=1$ corresponds to the 3rd order KdV-type equations; if $n_0=2,3$, respectively, the nonlinear evolution equations in the B\"acklund chart are all of the 5th and 7th order.

The  links among the corresponding members in the hierarchies can be depicted via the same B\"acklund 
chart in Section 5, that is, for each $n\in\N$, it holds 
\begin{gather*}
\mbox{$\footnotesize
\boxed{\! \left[\Phi_{1} (u)\right]^{n} \!u_x\!}\, {\buildrel (a)
\over{\text{\textendash\textendash}}} \, \boxed{\! \left[\Phi_{2} (v)\right]^{n}\! v_x\! } \, {\buildrel (b) \over{\text{\textendash\textendash}}} \, \boxed{\! \left[\Phi_3(w)\right]^{n}\! w_x\! } \, {\buildrel (c) \over{\text{\textendash\textendash}}} \,
 \boxed{\! \left[\Phi_{4} (\varphi)\right]^{n}\! \varphi_x\!}
 {\buildrel (d) \over{\text{\textendash\textendash}}}\, \boxed{\! \left[\Phi_{5} (s)\right]^{n}\! s_x \!} \,
{\buildrel (e) \over{\text{\textendash\textendash}}} \, \boxed{\!\left[\Phi_6 (\rho)\right]^{n-1} \rho^3\rho_{xxx}\!}\,
$}\label{BC3*}
\end{gather*}
where,  the recursion operators\footnote{The recursion operators here listed are well known ones with the 
only excepition of \eqref{psiw}, see, for instance \cite{CalogeroDegasperis}.}  are, respectively
{\begin{alignat*}{5} 
 &  \Phi_1(u)\equiv &  \Phi_{\text{KdV}} (u) & =D^2+2DuD^{ -1}+2u &    \text{(KdV), }&    \\
 &  \Phi_2(v)\equiv & \Phi_{\text{mKdV}} (v) & =D^2-4DvD^{ -1}vD &    \text{(mKdV),} &     \\
 &  \Phi_3(w)\equiv & \Psi (w) & =
{1 \over {2 w}}D w^2 \!\left[ D^2+2 U + D^{ -1} 2 UD\right]\! {1 \over {w^2 }} D^{ -1} \!2w &    \text{(new eq), }&    \\
 &  \Phi_4(\varphi)\equiv  & \Phi_{\text{KdVsing}} (\varphi) & = 
\varphi_x\left[ D^2 +  \{ \varphi ; x\}+D^{-1}  \{ \varphi ; x\}D\right] {1\over \varphi_x}  &  \text{(KdV sing.), }&    \\
 &  \Phi_5(s)\equiv & \Phi_{\text{KdVsol}} (s) & =D s\left[ D^2 +   S + D^{ -1}  SD\right] {1\over s} D^{-1} &    
 \text{(int.\ sol KdV),} &  
\\
 &  \Phi_6(\rho)\equiv & \Phi_{Dym} (\rho) & =\rho^3D^3 \rho D^{-1}\rho^{-2} &    \text{(Dym), }&   \end{alignat*}}
where 
\begin{equation}
U:=  {w_{xx} \over {w}}- 2 {w_{x}^2 \over {w^2}}~~~~,~~~~S:=  ( {s_x \over s } )_x -{1 \over 2}   ( {s_x \over s}) ^ 2~~.
\end{equation}
and the links among such  hierarchies of nonlinear evolution equations, are indicated in \eqref{links}.

\medskip
\noindent {\bf Remark} \\ \noindent 
{ The Dym hierarchy, generated on application of the hereditary recursion operator
\cite{CalogeroDegasperis, Fuchssteiner:Carillo:1989a, Lou, Leo} to the Dym equation,
\begin{equation}\label{dym}
 \rho_t=\left[ \Phi_{Dym} (\rho)\right]^n \rho^3\rho_{xxx}~, n\ge 0, ~~ \text{where}~~ \Phi_{Dym} (\rho)=\rho^3D^3 \rho D^{-1}\rho^{-2}
\end{equation}
is connected to all the hierarchies in 
which appears in the B\"acklund charts in Section \ref{ext-bc}; hence, it is also related to 
the hierarchy of nonlinear evolution equation whose base member is  \eqref{new}.
}

\section{Remarks, perspectives and open problems}
\label{rem}
This Section is devoted to collect some remarks and open problems which arise from the present  results.
This study is strictly connected with  the research program which involves C. Schiebold, together with, 
lately, also M. Lo Schiavo  and E. Porten, and the author concerning  operator evolution equations and 
their properties.  The results 
already obtained, based on the theoretical foundation in \cite{Carl Schiebold} and references therein,   
concern KdV-type non-Abelian equations in \cite{Carillo:Schiebold:JMP2009,  Carillo:Schiebold:JMP2011, SIGMA2016,  JMP2018} where, analogies as well as some  notable differences  which arise in the 
non-Abelian case are  pointed out. 
Non-Abelian Burgers hierarchies are studied  in  \cite{SIMAI2008, Carillo:Schiebold:JNMP2012, MATCOM2017}.
Comments concerning the comparison non-Abelian vs. Abelian results are  comprised in 
 \cite{Carillo:Schiebold:JMP2009, Carillo:Schiebold:JMP2011, ActaAM2012, 
SIGMA2016, MATCOM2017}. A remarkable property to stress in the present contest is that, also 
in the non-Abelian case, hereditariness is preserved under  B\"acklund  transformations. Hence, proved  the 
hereditariness of one  recursion operator \cite{Schiebold2010, Carillo:Schiebold:JNMP2012}, the
 hereditariness of all the recursion operator of other non-Abelian equations linked to it follows, see \cite{Schiebold-6dic1, Schiebold-6dic2, Schiebold-6dic3, Carillo:Schiebold:JMP2011, SIGMA2016}.
 \bigskip 
 
\noindent Some remarks follow.
\begin{itemize}\itemsep=0pt
\item The  B\"acklund chart in Section \ref{ext-bc} extends that one in   
\cite{Fuchssteiner:Carillo:1989a, walsan2}  finds its analogous B\"acklund chart in 
\cite{BS1, Rogers:Carillo:1987b}, where 5th order nonlinear evolution equations, i.e. 
Caudrey-Dodd-Gibbon-Sawata-Kotera (CDG-SK) \cite{CDG, SK}
and Kaup-Kupershmidt (KK) equations, which, in turn, play the role of the KdV and mKdV equations
appear. The  5th order nonlinear evolution equation analogous to the  Dym eq is
the Kawamoto equation \cite{Kawa}: it is linked via the reciprocal transformation 
\eqref{rec} to the singularity manifold equation \cite{Weiss} related to the CDG-SK equation. 
However, further to the many analogies the transformations  the two different B\"acklund charts 
are not exactly the same. Hence, 
the question arises whether or not there exist a 5th order analog of equation \eqref{new}. 
\item All the structural pro\-perties which are preserved under B\"acklund transformations are 
enjoyed by all nonlinear evolution equations in the same B\"acklund chart.
 This is the case, in particular, of hereditariness of recursion operators. Also the 
 Hamiltonian and/or bi-Hamiltonian structure \cite{[12], {Fuchssteiner:Carillo:1990a}, 
 Benno-Walter, Magri} in preserved under  B\"acklund transformations. Hence, the  Hamiltonian 
structure admitted by equation  \eqref{new} \cite{boris90} can be recovered 
from the presented B\"acklund chart.
\item In \cite{walsan1} the   $(2+1)$ dimensional KP, mKP and Dym hierarchies are connected via  B\"acklund
transformations; then in \cite{walsan2} it is shown that, when suitable constraints are imposed, the B\"acklund chart in \cite{Fuchssteiner:Carillo:1989a} is obtained. A further question which arises is whether the $2+1$ KP  eigenfunction equations \cite{boris90}   can be included in the B\"acklund chart constructed in \cite{walsan1}.  
\end{itemize}

Further questions concerning open problems and perspective investigations in the case of operator 
equations are  referred to \cite{SIGMA2016}. 

\subsection*{Acknowledgements}
The author wishes to thank with gratitude Boris Konopeltchenko for helpful discussions.

The financial support of G.N.F.M.-I.N.d.A.M.,  I.N.F.N. and \textsc{Sapienza}  University of Rome, 
Italy are gratefully acknowledged.

\end{document}